%% file: Haertig-3-submitted-Arxiv.tex
\documentclass[11pt]{article}
\usepackage{amsmath,amssymb,bbm,pslatex,xcolor,graphics}

\usepackage{xcolor}
\definecolor{blue}{rgb}{0,0,1}
 \usepackage[pagebackref,bookmarks]{hyperref} 
\hypersetup{pagebackref=true,bookmarks=true,
       hyperfigures=true,
	unicode=true,
	colorlinks=true,
	citecolor=blue,
	linkcolor=blue,
	anchorcolor=red
}
\usepackage[textsize=footnotesize,color=blue!40, bordercolor=white]{todonotes}

\usepackage{etoolbox}
\makeatletter
\def\new@mathgroup{\alloc@8\mathgroup\mathchardef\@cclvi}
\patchcmd{\document@select@group}{\sixt@@n}{\@cclvi}{}{}
\patchcmd{\select@group}{\sixt@@n}{\@cclvi}{}{}
\makeatother

\usepackage{lmodern}
\usepackage{mathpazo}
\usepackage{fourier} 

\usepackage{amsmath}
\usepackage{amssymb}

\usepackage{graphicx,bbm,color,calc,ifthen,theorem,capt-of,tocloft}
\DeclareMathSymbol{\gameb}{\mathord}{AMSb}{"61}

\usepackage{amsmath}
\usepackage{amssymb}
\usepackage{fancyhdr}

\include{pdwstyle}
\def\sset{\mbox{ $\subseteq$ }}
\catcode`\<=\active \def<{
\fontencoding{T1}\selectfont\symbol{60}\fontencoding{\encodingdefault}}
\catcode`\>=\active \def>{
\fontencoding{T1}\selectfont\symbol{62}\fontencoding{\encodingdefault}}
\catcode`\|=\active \def|{
\fontencoding{T1}\selectfont\symbol{124}\fontencoding{\encodingdefault}}

\newcommand{\tmop}[1]{\ensuremath{\operatorname{#1}}}

\newcommand{\tmmathbf}[1]{${\boldmath{$ #1$}}$}

\newtheorem{theorem}{Theorem}[section]

\newtheorem{corollary}[theorem]{Corollary}
\newtheorem{definition}[theorem]{Definition}

\newtheorem{lemma}[theorem]{Lemma}

\def\nod{\noindent}

\def\mb{\mbox{}}

\def\nod{\noindent}

\def\ed{\end{document}}
\def\mi{\mathsf{I}}

\def\bth{\begin{theorem}}
\def\eth{\end{theorem}}

\def\blem{\begin{lemma}}
\def\elem{\end{lemma}}

\def\bdefn{\begin{definition}}
\def\edefn{\end{definition}}

\renewcommand{\qed}{\mbox{ } \hfill Q.E.D.}

\def\mid{\, | \, }

\def\wt{\widetilde}

\def\m``{\mbox{``}}
\def\mq{\mbox{''}}

\def\pr{P\v r\'ikr\'y }
\def\ha{H\"artig }

\def\mp{\mathbbm{P}}
\def\mpn{\mathbbm{P}^{\nu}}

\def\co{\pa{h,H}}
\def\sp{\mathop{sp}}

\def\lep{L[E^{P}]}
\def\le0{L[E_{0}]}
\def\aln{\alpha < \nu}

\def\a{\alpha}
\def\k{\kappa}
\def\n{\nu}
\def\b{\beta}
\def\th{\theta}
\def\ii{\iota}

\def\l{\lambda}
\def\g{\gamma}

\newcommand\lv[2]{\lambda_{\omega {#1} + { #2} }}
\def\om{\omega}
\def\ups{\upsilon}
\def\forces{\Vdash}

\def\cs{\mbox{$C^{\ast}$}}
\def\ok{\mbox{$O^{k}$}}
\def\mbq{\mbox{``}}
\def\ki{\k_{\ii}}
\def\mi{M_{\ii}}
\def\io{\iota_{0}}
\def\ks{\mbox{$K^{\ast}$}}
\def\mm{{\em mutatis mutandis}}
\def\colr{\color{red}}
\def\power{{P}}

\begin{document}

\title{Closed  unbounded classes and the H\"artig quantifier model }
\author{ P.D. Welch\\  School of Mathematics, \\University of Bristol,\\
  Bristol, BS8 1TW, England}
\date{24.xii.18}

\maketitle

\begin{abstract}
 
 
We show  that assuming modest large cardinals, there is a definable class of ordinals, closed and unbounded beneath every uncountable cardinal, so that for any closed and unbounded subclasses  $P, Q$,   $\pa{L[P],\in ,P }$ and $\pa{L[Q],\in ,Q}$ possess the same reals, satisfy the Generalised Continuum Hypothesis, and moreover are elementarily equivalent. The theory of such models is thus invariant under set forcing. They also all have a rich structure satisfying many of the usual combinatorial principles and a definable wellorder of the reals. One outcome is that we can characterize the inner model constructed using definability in the language augmented by the H\"artig quantifier when such a $P$ is itself $Card$.
  \footnote{We should like to warmly thank the authors of \cite{KeMaVa2016} for many discussions on their paper. We first presented this result at the CIRM-Luminy meeting in Oct. 2017, and also should like to thank I. Neeman for pointing out an egregious and nonsensical  error in a version of the main theorem here claimed in our talk.}
\end{abstract}



\section{Introduction}

In this paper we consider inner models of the $ZFC$ axioms using constructibility relative to a predicate consisting of a closed and unbounded (cub) class of ordinals. Such models, so of the form $\pa{L[P],\in, P}$, are easily defined (see Kanamori \cite{Ka}). There are a number of questions one may ask about such: what structural properties they may have: are they models of $GCH$? of $V=HOD$? Does $\square$ hold in them? How do they relate to other well known inner models - are they fine structural? What are their reals? What are their grounds?

Of course if the universe is too thin, these dissolve into triviality, for example if $V=L$ in the first place. Forcing constructions over $L$ also give some not terribly interesting consistency results. However it turns out that with a modest large cardinal assumptions in the universe (that there is a measurable limit of measurable cardinals, or more precisely  that there exists an elementary embedding of an inner model with a proper class of measurable cardinals to itself - we'll call the latter assumption $O^{k}$ ($= O^{kukri}$))  then we have the following perhaps surprising theorem:

\bth \label{Thm1.2} $(ZFC)$ Suppose $O^{k}$ exists.
There is a definable proper class $C\sset On$ that is cub beneath every uncountable cardinal, so that for any definable cub subclasses $P,Q \sset C$:

$$\mathbbm{R}^{L[P]}=\mathbbm{R}^{L[Q]}; \quad \pa {L[P],\in, P} \equiv \pa {L[Q],\in, Q}
$$ 
where the elementary equivalence is in the language $\call_{\dot \in, \dot P}$ with a predicate symbol $\dot P$.  Consequently these models are all similar to one another: they have the same reals, and their theories are invariant under set forcing.
\eth

One might {\em prima facie} have surmised that a clever choice of elements of $P$ might have allowed some coding of interesting sets in order that at the very least the theories of two such models would be different. But apparently not. A particular example of course is when $P=Card$ itself, the latter the class of uncountable cardinals. These models $L[P]$ all have a rich structure and we have a complete picture of them: they can be considered as a form of generalised \pr class generic extensions of a fine structural  model with a proper class of measurables (hence the need for the hypothesis $O^{k}$). This fine structural model will naturally  form the {\em core model} of the class $L[P]$, for such $P$. They thus have nice  combinatorics: $\square_{\lambda}$ holds everywhere, $GCH$ holds. They all have the same set of reals. The elements of $P$ are all J\'onsson in the model $L[P]$, but not much more. (See Cor. \ref{Cor2.4} below for a listing of such properties.) We should point out that these results can be extended easily to considering sequences $P,Q$ from $C$ of the same bounded but limit order type order type: again the displayed formulae of Theorem 1.1 would hold for such $L[P]$ and $L[Q]$ too.

We apply this to solve the question of the identity of the H\"artig quantifier model (which was the starting point for this paper). In \cite{KeMaVa2016} the authors consider the possibilities of using the G\"odel method of defining a hierarchy of constructible sets, but augmenting the logic from straightforward first order definability to one where a new quantifier $\mathsf{Q}$ is added to the language. If the enhanced language is called $\call^{\ast}$ they build a model as follows:

$$
\begin{array}{lcr} L_0^{\mathsf{Q}}& = & \emp\\
L_{\alpha+1}^{\mathsf{Q}} &=& Def _{\call^{*}}(L_{\alpha}^{\mathsf{Q}}) \\
L_{\lambda}^{\mathsf{Q}} & = & \bigcup_{\alpha < \lambda}L_{\alpha}^{\mathsf{Q}}
\end{array} 
$$
and then $L^{\mathsf{Q}}=  \bigcup_{\alpha \in On}L_{\alpha}^{\mathsf{Q}}$.

If the quantifier $\mathsf{Q}$ is the \ha quantifier $I$, they dub the model $C(I)$. 

\bdefn The \ha quantifier $I$ has the following interpretation:
$$I x y \varphi(x,\vec p)\psi(x,\vec p) \equi |\{a : \varphi(a,\vec p) \}| =  |\{b : \psi(b,\vec p) \}|.$$
\edefn

For a summary  of facts concerning this quantifier see \cite{HKPV91}. It is an important point to note that the construction of an $L^{\mathsf{Q}}$-hierarchy in such cases feeds in information from $V$. We would not expect the construction of such a hierarchy to be in any way absolute. Other than in trivial cases (such as when $V=L$) we should not expect that $(V=L^{\mathsf{Q}})^{L^{\mathsf{Q}}}$ for example.

The paper \cite{KeMaVa2016} shows, {\em inter alia}, the following results:\\

\nod$\bu$ If $L^{\mu}$, the least inner model of a measurable cardinal, exists then $L^{\mu} \sset C(I)$.\\
\bu \,$Con(ZFC + \exists \kappa (\kappa $ supercompact) $) \Rightarrow Con(ZFC + C(I)\neq HOD).$\\

However it is left open as to exactly what model $C(I)$ is, or of what kind.  It is easy to see that with $\mathsf{Q}=I$ that $C(I)=L[Card]$ where $Card$ is a one place predicate true of the infinite cardinals (and $L[Card]$ is defined from the usual first order relativised constructibility hierarchy from the predicate). But that alone tells one very little about the structure of $C(I)$ for example whether it has large cardinals, or the $(G)CH$ holds there. However $Card$ is a cub class contained in the $C$ of the theorem above. The theorem and its proof are thus applicable to $L[Card]$. It is important to require the closure of the classes $P$ of Theorem \ref{Thm1.2}: let ${Reg}$ be the class of regular cardinals; then we can show that $L[Card]$ is, barring trivialities, a proper subclass of $L[Reg]$.

We shall have:
\setcounter{section}{5}

\setcounter{theorem}{0}
\begin{theorem}\label{KinC(I)} 
 $\neg O^{k} \Equi K^{I} = K$.
\end{theorem}

\nod where we shall set $K^{I}= (K)^{C(I)}$. Here $K$ is the {\em core model}, which we regard as here constructed {\em \`a la Jensen} for which see the original manuscript of \cite{Je89}, where the discussion is about mice with measures of order zero, which is all that we shall deal with here. Similarly the first part of \cite{Z02} gives a full exposition of this theory.  Such a model is one of a family of models of the form $\pa{L[E],\in,E}$ where $E$ is a {\em coherent sequence of extenders}. In this context the extenders can be rendered as simply filters (again see \cite{Z02}). These models possess {\em fine structure}, have global wellorders of their domains, satisfy a strong form of the $GCH$ and have strong combinatorial properties, such as Jensen's $\Box_{\k}$-property everywhere. For `small' or `thin' $L[E]$ models, they will, like $L$, be models of the statement `I am $C(I)$':

\begin{corollary} $(V=L[E])$ \,\, $\neg O^{k}\equi V=C(I)$.
\end{corollary}

\setcounter{section}{1}

\setcounter{theorem}{1}
This note, assuming large cardinals, rather just that $O^{k}$ exists, identifies this model: $C(I)$ is a generalised \pr forcing extension of (an iterate of) the smallest inner model with a proper class of measurable cardinals. One way to express this is to say that, for limit ordinals $\lambda$ the $\omega$-sequences of successor cardinals $c_{\lambda }\dfs \pa {\aleph_{\lambda +i}\,\mid\, 0<i<\omega}$ form \pr -sequences for the $L[E]$ model which is the least inner model with a measurable cardinal on every $\aleph_{\lambda + \omega}$.  We do this in such a manner that the class $\pa{c_{\lambda}\mid \lambda \in Lim}$ is $\mathbbm{P}^{\infty}$-generic over the model $L[E]$ for a certain class forcing $\mathbbm{P}^{\infty}= \mathbbm{P}^{Card,\infty}$. The source of this forcing is Magidor's iterated \pr-forcing (\cite{Ma74} or see \cite{Gi07}) which has a full support; however as the measures in the model $L[E]$ are sparsely distributed (there are inaccessible limits of measurables, but no measurable limits of measurables) the forcing can simplified. Here we use such a simplified version as was used for countable sequences in \cite{W2}, but more relevant here, for any set sized sequences of measurables - again with no measurable limits of measurables - analysed in detail by Fuchs \cite{Fu05}. That $C(I)\neq HOD$ will now follow from the existence of $O^{k}$ (Cor. 5.5).\\

In a final section we make a few remarks about the relationship between $C(I)$ and $C^{\ast}$ - the latter the inner model defined using the additional `cofinality $\om$'-quantifer $\mathsf{Q}^{\om}$. ($\cs$ is coextensive with $L[\tmop{Cof}_{\om}]$ where $\tmop{Cof}_{\om}$ is  the class of ordinals of cofinality $\om$.)
There is extensive discussion in \cite{KeMaVa2016} on this model. A model may be large in one sense, even if it does not have any, say measurable cardinals, of its own: it may have inner models with very large cardinals instead, and this would surely count as the model being `large'. However in all of the results there, some of which assume very large cardinals in $V$, the outcomes for $\cs$ are nevertheless all consistent with it being also a thin model. We show here that it must be larger than $C(I)$, but not by much, only in that $O^{k}\in C(I)$. So, one might conjecture that $\cs$ is also thin:\\

{\em Conjecture}: $\cs$ does not contain a mouse with a measurable of Mitchell order $\om_{1}$. Or alternatively no mouse with a measurable limit of measurable cardinals with Mitchell order $\omega_{1}$.\\
 
\nod Our result here does not imply that a mouse with a measure of Mitchell order $> 0$ is in $C(I)$.

\section{The model $L[E_{0}]$}

The principal model $L[E_{0}]$ we shall use can be derived as follows.
\bdefn
Let $O^{k}$ name\footnote{``O Kukri'' - from a Ghurka weapon somewhat intermediate between a dagger and a sword.}  $M_{0}$ being the least sound active mouse of the form $
M_{0}\dfs \pa { J^{E^{M_{0}}}_{\alpha_{0}}, E^{M_{0}}, F_{0}}$ so that
$$ M_{0}\models \mb`` F_{0} \mbox{ is a normal measure on $\kappa_{0}$ } \wedge \forall \tau <\kappa_{0}\exists \lambda < \kappa_{0}( \tau < \lambda 
 \mbox{ and $\lambda$ a measurable cardinal } )
\mbox{''.}$$
\edefn 

Here we mean a mouse in the sense of \eg \cite{Z02}, and the $E^{M_{0}}$ sequence is a coherent sequence of filters from which we are constructing. Then the following list of {\em Facts} are either common knowledge or are easily derived from standard arguments:\\
(i) $M_{0}$ is a countable mouse with $\rho^{1}_{M_{0}} = \omega$ - the first projectum drops to $\omega$ and there is a $\Sigma_{1}^{M_{0}}$ definable map of $\omega$ onto $J^{E_{M_{0}}}_{\alpha_{0}}$.\\
(ii) We may form iterated ultrapowers of  $M_{0}$ repeatedly using the top measure $F_{0}$ and its images to form iterates $
M_{\iota}\dfs \pa { J^{E_{M_{\iota}}}_{\alpha_{\iota}}, E_{M_{\iota}}, F_{\iota}}$ so that
$ M_{\iota}\models$ ``$ F_{\iota} \mbox{ is a normal measure on }\kappa_{\iota}$''.\\
(iii) These iterations generate, or ``leave behind'', an inner model $L[E_{0}] \dfs \bigcup _{\iota \in On } H_{\kappa_{\iota}}^{M_{\iota}}.$\\
(iv) The cub class of critical points $C_{M_{0}}=\pa{ \kappa_{\iota}\mid \iota \in On} $ forms a class of indiscernibles that is cub beneath each uncountable cardinal, for the inner model $L[E_{0}]$. Indeed an elementary skolem hull argument shows that the  $\{\kappa_{\iota}\}$ form a class of {\em generating indiscernibles} for $L[E_{0}]$ just as the Silver indiscernibles from $O^{\sharp}$ do for $L$.\\
(v) From (iii) we have that for any $\iota <\nu$ that $L_{\kappa_{\iota}}[E_{0}] \prec 
L_{\kappa_{\nu}}[E_{0}] \prec L[E_{0}]$. Moreover for any $\iota \in On$ we have that $H_{\kappa_{\iota}^{+}}^{L[E_{0}]}= 
|J^{E^{M_{\iota}}}_{\alpha_{\iota}}|
$, where $\kappa_{\iota}^{+}$ is the successor cardinal of $\kappa_{\iota}$ in the sense of $L[E_{0}]$ and is thus identical to $\alpha_{\iota}$. If $j_{\iota,\nu}: M_{\iota} \longrightarrow M_{\nu}$ is the iteration map between the iterates displayed, we shall thus have that also $j_{\iota,\nu}:  H_{\kappa_{\iota}^{+}}^{L[E_{0}]}\longrightarrow H_{\kappa_{\nu}^{+}}^{L[E_{0}]}$ is an elementary embedding, which extends to an elementary map
$\widetilde \jmath _{\iota,\nu}:  {L[E_{0}]}\longrightarrow {L[E_{0}]}$.
 (Again this is similar to the corresponding fact in the embeddings of $L$ coming from iterations of the ``$O^{\sharp}$-mouse'': for $\mu< \nu$ Silver indiscernibles for $L$, we have an elementary map $j_{\mu,\nu}: L_{\mu^{+}}\longrightarrow L_{\nu^{+}}$, which extends to a map $\widetilde \jmath_{\mu,\nu}: L\longrightarrow L$.)\\
(vi) We may if we wish think of $L[E_{0}]$ to have the same domain as the model $L[\vec U]$ where $\vec U$ is a sequence of filters on the $\kappa_{\iota}$ which are normal measures in $L[\vec U]$. The fine structure of the latter model was originally done {\em \` a la} Dodd-Jensen (\cite{D}) rather than the $L[E]$ style of Jensen in \cite{Z02}. But the models have the same domain of sets.

We call a class $P$ of ordinals {\em appropriate} if $P\sset C_{M_{0}}$ is  closed and unbounded.
For such an appropriate $P$ let $ \pa {\lambda_{\iota} | \iota \in On }$ be its strictly increasing enumeration. 
 Further, for $\alpha\in On$ we set $c({\alpha}) = c^{P}({\alpha}) =\pa {\lambda_{\omega\alpha +k}\mid 0<k<\omega}$  and $\mu _{\alpha}=\mu _{\alpha}^{P} \dfs \lambda_{\omega \alpha+\omega}$.  Note the particular case of interest for later is the appropriate class $P=Card$. 
 With this notation then we shall see the following:

\bth Assume that $O^{k}$ exists and $P$ is an appropriate class.
 (i) $K^{L[P]}= L[ E^{P}]$ where $ E^{P}$ is a coherent filter sequence so that $L[ E^{P}]\models $`` $\kappa$ is measurable'' if and only if $\kappa = \mu_{\alpha}$ for some $\alpha$. \\ 
 (ii)
The class $ \pa {c^{P}({\alpha})}$ $\dfs\pa{c^{P}({\alpha})\mid \alpha\in On}$ of $\omega$-sequences is mutually \pr-generic over $L[E^{P}]$ for the forcing $\mathbbm{P}^{P} 
$ and $L[P] = L[E^{P}][ \pa {c^{P}({\alpha})}]$.
\eth

A corollary of (the proof of) our theorem will be the following (a restatement if Theorem 1.1):
\bth Assume that $O^{k}$ exists.
Let $P,Q\sset C_{M_{0}}$ be any two appropriate classes. Then
$$\mathbbm{R}^{L[P]}=\mathbbm{R}^{L[Q]}; \quad \pa {L[P],\in, P} \equiv \pa {L[Q],\in, Q}
$$ 
where the elementary equivalence is in the language $\call_{\dot \in, \dot P}$ with a predicate symbol $\dot P$.
\eth

\begin{corollary}\label{Cor2.4} Assume $O^{k}$ exists.
Let $P$ be any appropriate class. Then in $L[P]$:\\
(i) Each $\mu_{\alpha}$ is J\'onsson,  and ${c_{\alpha}}$ forms a {\em coherent  sequence of  Ramsey cardinals} below $\mu_{\alpha}$ (see Koepke {\em\cite{K84}}). But there are no measurable cardinals.\\
(ii) For any $L[P]$-cardinal $\kappa$ we have $\Diamond _\kappa$, ${\Box}_\kappa$, $(\kappa,1)$-morasses \etc \etc\\
(iii) The $GCH$ holds but $V\neq HOD$.\\
(iv) There is a $\Delta^{1}_{3}$ wellorder of $\mathbbm{R}= \mathbbm{R}^{K^{L[P]}}$; $Det(\alpha$-$\tmmathbf{\Pi^{1}_{1}})$ holds for any countable $\alpha$ (see {\em \cite{W2}}), but $
Det(\tmmathbf{\Sigma^{0}_{1}}(\tmmathbf{\Pi^{1}_{1}}))$ fails (Simms, Steel, see {\em\cite{St82}}). 
\end{corollary}
Indeed anything else that holds after a \pr \!-generic extension of the $L[E^{P}]$ model. Notice that $(Card)^{L[Card]}$ will be very far from $Card$ as any $\mu \in Card$ will be in $L[Card]$ a Ramsey cardinal (hence weakly compact) or a limit of such.\\

We note the following for later use.
\blem\label{minimality}
Suppose $\ok$ exists. Let $L[E_{0}]$ be the model defined above. Let $L[E']$ be any other model with a proper class of measurable cardinals, with $L[E_{0}]=^{\ast}L[E']$ in the mouse/weasel ordering. Then $L[E']$ is a simple iterate of $L[E_{0}]$.
\elem
\pf As the models are $=^{*}$ equivalent the comparison of  the models will be simple iterations on both sides. The claim is that the iteration on the $L[E']$ side is trivial, \ie no ultrapower is ever taken. However note that if $N_{0}$ is the least sound mouse that generates $L[E']$ then $N_{0}= M_{0}= \ok$.
\qed\\

In one obvious sense then $L[E_{0}]$ is the `minimal' model with a proper class of measurable cardinals.

Woodin in \cite{Wo96} considered the question of what occurs when an $\omega$-sequence of ordinals is added to $L$. Besides reals added by forcing of course, much can happen. He shows that if $s$ is an $\omega$-sequence of ordinals then $L[s]$ is a model of GCH. This also used a \pr-forcing and a short core model analysis. We comment below on what happens when we choose an $\omega$ sequence or indeed any limit length sequences $P$ contained in $C$.

\subsection{Universal Iterations}

We place here a general discussion on {\em universal iterations} of a mouse, which we shall use only here as a device to ensure that certain iterations of a model, although defined externally to the model, leave inaccessibles of the model fixed. These results appeared in a somewhat more {\em recherch\'e} form in \cite{SW}. 

\bdefn \label{defuniv}(\cite{SW} Def. 2.8)
Let $M$ be a mouse and $ \theta > On^{M}$ be an $M$-admissible ordinal. Then
$\pa{M_{\eta}, \pa{\ups_{\eta,\iota}}_{\eta \leq \ii \leq\theta}, {\kappa_{\eta}}}_{\eta \leq\theta}$ with indices $\pa{\nu_{\eta}}_{\eta \leq \theta}$
is an {\em $n$-universal iteration} of $M= M_{0}$ {\em of length } $\theta$, if (i) there are no truncations and dropping of degree of the iteration at any stage $\alpha < \theta$ and (ii) for any measure $F= E^{M_{\alpha}}_{\nu_{\alpha}}$ with $crit(F)< \rho^{n}_{M_{\alpha}}$ there is $\beta<\theta$, $\alpha < \beta$ with $E^{M_{\beta}}_{\nu_{\beta}}= \ups_{\alpha,\beta}(F)$.

\edefn
Thus, in an universal iteration, every extender  (or its image under the iteration so far) appearing on any extender sequence of the iteration is used unboundedly often before $\theta$. We shall be using the simplified version of the above where $n=\omega$ and the extenders are measures are all elements of the models appearing, which are themselves $ZF^{-}$ models (and so $ \rho^{\om}_{M_{\alpha}}= On \cap M_{\alpha}$ throughout). The next lemma states that, although there can be many universal iterations of given length, any two such end up with isomorphic results.

\bth (\cite{SW} Thm. 2.9)
Let $\theta> On^{M}$ be an $M$-admissible ordinal. If  \,
$\calu = \pa{M_{\eta}, \pa{\ups_{\eta,\iota}}_{\eta \leq \ii \leq\theta}, {\kappa_{\eta}}}_{\eta \leq\theta}$ and  $\,\overline{\calu} = \pa{\bar M_{\eta}, \pa{\bar\ups_{\eta,\iota}}_{\eta \leq \ii \leq\theta}, {\bar \kappa_{\eta}}}_{\eta \leq\theta}$ are any two $n$-universal iterations of $M=M_{0}=\bar M_{0}$
of length $\theta$
then $M_{\theta} = \bar M_{\theta}$. 
\eth

We may define a universal iteration in $L[M]$:

\blem (\cite{SW})
Let $\theta< \th_{0}$ be two $M$-admissible ordinals. Then there is an $n$-universal iteration of $M$ up to $\th$ which is an element of $L_{\th_{0}}[M]$.
\elem

The point of a universal iteration is that any other iteration of the first model of a shorter length can be embedded into the universal iteration. We formulate that as follows.

\bth (\cite{SW} Thm. 2.10)\label{univ}
Let $\theta$ be an $M$-admissible ordinal. If \,
$\calu = \pa{M_{\eta}, \pa{\upsilon_{\eta,\iota}}_{\eta \leq \ii \leq\theta}, {\kappa_{\eta}}}_{\eta \leq\theta}$
 is an $n$-universal iteration of $M=M_{0}$ up to $\theta$,
and $ \calj = \pa{ {N_{\ii}}, \pa{\pi_{\ii ,j}}_{\ii\leq j\leq\mu}, {\kappa_{\ii}} }_{\ii \leq \mu}  $ is any length $\mu+1<
\th$ $n$-iteration of $M= N_{0}$, (with no truncations or drops in degree) then there is an iteration
$\cali = \pa{ {P_{\ii}}, \pa{\sigma_{\ii ,j}}_{\ii\leq j\leq\xi}, {\kappa_{\ii}} }_{\ii \leq \xi}  $ of  $P_{0}= N_{\mu}$ of length some $\xi+1< \th$ (with no truncations or drops in degree) so that for some $\beta$, $P_{\xi}= M_{\beta}$.
\eth

We thus say that a universal iteration of length $\theta$ absorbs all shorter length (appropriate) iterations of the first model.
 We shall only use this construction in one particular case. Let $N$ be an inner model with only boundedly many measurable cardinals, bounded by some least $N$-regular cardinal $\theta$ say. Then we may just as easily as above define a $\theta+1$ length universal iteration of the proper class $N$ using the measures which are all below $\theta$, and moreover we define this universal iteration {\em in} $N$. But to make it about sets, we consider just some sufficiently large initial segment $N\rest\gamma$ where $\gamma$ is an $N$-inaccessible limit of $N$-inaccessibles. (Our intended $N$ will satisfy there are such.)  We thus consider the universal iteration to be on the first model $Q_{0}= N\rest\gamma$ of the universal iteration 
 $ \pa{Q_{\eta}, \pa{\upsilon_{\eta,\iota}}_{\eta \leq \ii \leq\theta}, {\kappa_{\eta}}}_{\eta \leq\theta}$

We then shall have: 

\blem\label{fixed1.5}
Let $\pi: N_{0}\longrightarrow N_{\io}$ be any simple iteration of $N_{0}=N\rest \gamma$ with $\gamma$ as above, of length $\io+1$. Then for any $N$-inaccessible $\bar\gamma \in (\theta,\gamma)$, $\bar\gamma$ is a fixed point of $\pi$: $\pi(\bar\gamma)=\bar\gamma$.
\elem
\pf As $Q_{0}$ has inaccessible height in $N$, $ \ups_{0,\theta}\m``\gamma \sset \gamma$ and indeed $\bar\gamma= \ups_{0,\theta}(\bar\gamma)$ for any $N$-inaccessible $\bar\gamma$ in our chosen interval. (Proven by induction on $j\leq \theta$ for the maps $\ups_{0,j}$
by the usual counting of functions in the internally defined iteration $\calu$.)
Further by the Theorem \ref{univ} there is an iteration $\cali  = \pa{ {P_{\ii}}, \pa{\sigma_{\ii ,j}}_{\ii\leq j\leq\xi}, {\kappa_{\ii}}  }_{\ii \leq \xi} $ of  $P_{0}= N_{\io}$ of length some $\xi+1< \th$ 
so that for some $\beta< \theta$, $P_{\xi}= Q_{\beta}$. However we have commuting maps $\ups_{0,\beta} = \sigma_{0,\xi}\circ \pi: Q_{0}\longrightarrow Q_{\beta}$.  But $\ups_{0,\beta}(\bar\gamma)=\bar\gamma$ as the $N$-inaccessibles are fixed points of these maps defined in $N$. So then $\pi(\bar\gamma)=\bar\gamma$ too.
\qed

\section{The Generalized {\pr} forcing}
In \cite{Fu05} is developed a style of iterated \pr  forcing intended for use when there are no measurable limits of measurables. This considerably simplifies the format of the forcing as the man\oe uvres needed for names in the full Magidor iteration of \cite{Ma74} are not needed. Moreover Fuchs proves a Mathias like characterisation (see Thm. \ref{mathias} below) which we shall make use of. The subsection 3.2 thus borrows heavily from \cite{Fu05}, but we shall adopt notation appropriate for this specific case.

\subsection{The model $L[E^{P}]$}
We first defined a simple iteration of $M_{0}$ by its top measure and its images used $On$ times, that left behind the inner model $L[E_{0}]$.  We fix an appropriate class $P$ for this whole discussion. We may then define a normal iteration of $L[E_{0}] \longrightarrow \lep$ to line up the measures of $L[E_{0}]$ onto the simple limit points of $P$, the $\pa{\mu^{P}_{\a}}_{\a \in On}$. We can reorganise these two into a single normal iteration.
$\pa{M_{\iota}, \pa{\pi_{\eta,\iota}}_{\eta \leq \iota<\infty}, {\kappa_{\iota}}}_{\iota\in On}$ where as usual at limit stages direct limits are taken.
 Indeed given the model $\lep$, the comparison coiteration of $(M_{0},\lep)$ (see \cite{Z02})  tells us what that iteration is by simply observing the $M_{0}$-side, as the $\lep$ model does not move in this. This iteration of $M_{0}$ `leaves behind' $\lep$. Between ultrapowers where the top filter from the relevant model is used are the intermediate ultrapowers lining up each of the full measures with the appropriate $\mu_{\a}$.  It is useful to identify the stages where the top measure is used: we let $C= C_{P} $ be this class of indices. It is easy to see that $C\sset C_{M_{0}}$ and is also cub in $On$.  Thus with $\eta <\iota$ both in $C$ we shall have Fact (v)  above (and the comments following) holding in this context \ie
we have that for any $\iota <\nu$ both in $C$, with $\pi_{\iota,\nu}: M_{\iota} \longrightarrow M_{\nu}$: \\

(1)
 There is an extension of $\pi_{\ii,\nu}$ to $\widetilde \pi_{\ii,\nu}$ with  $\widetilde \pi_{\ii,\nu}: \lep \longrightarrow_{e} \lep
$.\\

Consequently we also have the $\pa{\kappa_{\ii}}_{\ii\in C} $, which are $\Sigma_{1}$-indiscernibles for the $M_{\iota}$, will be full indiscernibles for $\lep$, and {\em inter alia} that

$$L_{\kappa_{\iota}}[E^{P}] \prec 
L_{\kappa_{\eta}}[E^{P}] \prec L[E^{P}].$$ 

We shall thus have that also $\pi_{\iota,\nu}:  H_{\kappa_{\iota}^{+}}^{L[E^{P}]}\longrightarrow H_{\kappa_{\nu}^{+}}^{L[E^{P}]}$ is a fully elementary embedding by noting that the domain of $M_{\ii}$ is precisely  this $H_{\kappa_{\iota}^{+}}^{L[E^{P}]}$ in the model being left behind.
We have then  that for each $\nu \in C$ that it is an inaccessible limit of measurables in $\lep$.

From the above, in $L[E^{P}]$ we have that $\pa{\mu_{\alpha} \mid \alpha \in On}$ is a proper class of discrete measurable cardinals with normal measures $U_{\alpha}$ (which are indexed on the $E^{P}$-sequence by $(\mu_{\alpha}^{+})^{L[E^{P}]}$ although that is not of much consequence in what follows).  
We note also the following:
\blem Fix $\kappa \in C_{M_{0}}$.
Let $\cali = \pa {N_{\ii}, \pa{\sigma_{\ii,}}_{\ii\leq\theta} }$ where $\theta < \kappa$ be a simple iteration of $N_{0}= L[E_{0}]$. Then $\sigma_{0,\theta}(\kappa)= \kappa$.
\elem
\pf Firstly note $\k$ is strongly inaccessible in $\le0$ as it is indiscernible there. The iteration $\cali$ is divided into two parts: those measures used below $\k$ and those above. It suffices to note that if the iteration below $\k $ does not move $\k$ the rest of the iteration using critical points $\kappa_{k}\geq \k$ will not move $\k$ as, in particular, $\k$ is not measurable in in $\le0
$. So it suffices to consider only those $\cali$ with measures used below $\k$. However for such an iteration, although not necessarily internally definable in $\le0$, one shows by induction on $\theta$
that $\sigma_{0,\theta}$ cannot move $\k$ as $\theta <\k$ (\cf the arguments using universal iterations  in Lemma \ref{fixed1.5}).
\qed\\

As a consequence we have:

(2) Any $\kappa \in C_{M_{0}}$, is only moved in an iteration $\sigma_{0,\th}:L[E_{0}] \longrightarrow L[E^{P}]$ if $\th \geq \k$ and  for some $\l <\k$ we have $\sigma_{0,\th}(\l ) \geq \k$. \\

\subsection{The forcing}
We proceed to define the forcing in $L[E^{P}]$ up to the $L[E^{P}]$-inaccessible cardinal $\nu \in C$.
\begin{definition} For $\nu \in C$ let $\mathbbm{P}^{\nu}= \mathbbm{P}^{P,\nu}$ be the following set of function pairs $\pa{ h,H}$ so that:\\
(i) $H\in \prod_{\alpha< \nu} U_{\alpha},\, \dom(h)=\nu$ and $\sp(h)$ is finite, where the latter, the support of $h$, is: $\sp(h)\dfs \{ x\in \dom(h) \mid h(x)\neq \emp\}$.\\
(ii) $\forall \alpha \in \sp(h)\, h(i) \in [\mu_{\alpha}]^{<\omega}$,\\
(iii)  $\forall \alpha \in \sp(h)\, h(\alpha)\sset \min H(\alpha)$.\\
(iv)  $\forall \alpha \in \sp(h)\, \forall \beta <\alpha ( \mu_{\beta}< \min(h({\alpha}))$. \\

For $\pa {f,F}, \pa{g,G}\in \mp^{\nu}$ set $\pa {f,F}\leq \pa{g,G}$ iff  $\forall \alpha <\nu (\, f(\alpha)\supseteq g(\alpha)\,\wedge\, f(\alpha)\backslash g(\alpha)\sset G(\alpha))$.\\

\end{definition}
The reader will recognise that we are using a form of \pr forcing with full support up to $\nu$. (Those familiar with \cite{Fu05} will see that we have further simplified by only seeking \pr sequences of length $\omega$ in the generic extension.)
We have the following basic properties (3)-(7) from Fuchs \cite{Fu05} p.939.\\
{\em Facts}:\\
(3) For any $\pa{h,H}\in \mpn$, any $\alpha <\nu$, there is $\pa{h',H'} \leq \pa{h,H} \wedge |h'(\alpha)|>n$.\\

\nod For the remainder of these Facts we let
 $G^{\nu}$ be $\mpn$-generic over $ \lep$, and we define $c= c_{G^{\nu}}$ by $c(\alpha)=\bigcup\{h(\alpha)\mid \ex H \co \in G^{\nu}\}$ for all $\alpha < \nu$. \\
 
\nod (4) Then $c\in \prod_{\alpha <\nu}(\mu_{\alpha}\back \bigcup_{\beta < \alpha}\mu_{\beta})^{\omega}$.\\
 
\nod(5) 
  $G^{\nu} = G_{c} $ where the latter is $ \{ \co \in \mp\mid \all \alpha < \nu (h(\alpha)\mbox{ \em is an initial segment of } c(\alpha)\wedge c(\alpha)\back h(\alpha)\sset H(\alpha))\}.
  $\\
  
   The last then yields that $\lep[c]=\lep[G^{\nu}]$.\\
   
\nod  (6) $\mpn$ has the $\nu^{+}$- c.c. (and this is best possible).\\

 \nod (7) For every $X\in (\prod_{\alpha <\nu}U_{\alpha})\cap \lep$, the set $\bigcup_{\alpha <\nu} (c(\alpha)\back X(\a))$ is finite.\\
 
 We have the following crucial Mathias-like characterization of this product of forcings, stated in our terms:

\begin{theorem}[Fuchs \cite{Fu05} Thm. 1]\label{mathias}
A function $d\in  \prod_{\alpha <\nu}(\mu_{\alpha}\back \bigcup_{\beta < \alpha}\mu_{\beta})^{\omega}$ is $\mpn$-generic over $L[E^{P}]$ if and only if  for every $X\in (\prod_{\alpha <\nu}U_{\alpha})\cap \lep$,
 $\bigcup_{\alpha <\nu} (d(\alpha)\back X(\a))$ is finite.
\end{theorem}
The combinatorics of this argument are somewhat involved so we don't repeat this here. But a corollary to this, also observed by Fuchs, allows a version of weak homogeneity which we shall exploit later. Since we have full products but only finite supports and thus only finitely many \pr-stems, if $c$ is any $\mpn$-generic and $p$ any condition, there is a finite perturbation $d$ of $c$ with $p\in G_d$.  Using the Mathias characterisation of \ref{mathias} this is the idea behind

\begin{corollary}[\cite{Fu05} Cor.1]
Let $c$ be $\mpn$-generic over $\lep$. Let $p\in\mpn$. Then there exists a sequence $d$ which is $\mpn$-generic over $\lep$ so that:\\
(i) $| \bigcup_{\aln}(c(\alpha)\triangle d(\alpha))|<\omega $\, ;\\
(ii) $p\in G_{d}$.
\end{corollary}

But such a $d$ is in $\lep[c]$ and we have then this model equals $\lep[d]$. Consequently we have also:

\begin{corollary}
If $\varphi(v_{0}, \ldots, v_{n-1})$ is any formula and $\check a_{1},\ldots \check a_{n-1}$ any forcing names for elements of $\lep$ and $\dot\Gamma$ a name for $c_{G^{\nu}}$, and $p\in \mpn$ we have $$p\forces_{\mpn} \varphi(\dot\Gamma, \check a_{1}, \ldots, \check a_{n-1}) \Imp  \mathbbm{1}\forces_{\mpn} \varphi(\dot\Gamma, \check a_{1}, \ldots, \check a_{n-1}).$$

\end{corollary}

Again from Fuchs we have (8)-(9):

 \nod (8) For $\gamma <\nu$ (not necessarily in $C$) $\mpn$ can be decomposed as a product $\mp^{\gamma} \times \mpn_{\gamma}$ with elements of $\mp^{\gamma} $ functions with domain $\gamma$ and those in $ \mpn_{\gamma}$ with domain $[\gamma, \nu)$.\\

 \nod (9) Forcing with $\mpn$ preserves all cardinals and cofinalities excepting the measurable cardinals, which are made cofinal with $\omega$ by the addition of the generic function $c$.\\

We also have:

\nod (10) 
 Let $\sigma$ be a sentence of the forcing language and $p\in \mpn$ be a condition. Then there is a `pure' or `direct' extension $q\in \mpn$ with $q || \sigma$, $q$ deciding $\sigma$. That is if $p=\pa{h,H}$, then such a  $q$ is the form $\pa{h,H'}$ where $H'(\beta)\sset H(\beta)$ for all $\beta < \nu$. (See, Gitik \cite{Gi07} Lemma 6.2). Further $ \mpn_{\gamma}$ adds no bounded subsets of $\k_\gamma$ - the $\gamma$'th measurable cardinal of $L[E^{P}]$ ({\em ibid.} Sect. 6.)

\section{The class version: the full forcing $\mp^{\infty}=\mp^{P,\infty}$}







We may consider the forcings $\mp^{\ii},\mpn$ as above, as defined for such  $\ii,\nu\in C$, $\ii<\nu$, within $\lep$.\\

(11) $\widetilde \pi_{\ii,\nu} (\mp^{\ii}) = \mpn$.\\


(12) $c$ is $\mp^{\ii}$-generic over $\lep$ if it is so over  $H_{\kappa_{\iota}^{+}}^{L[E^{P}]}$.

\pf Let $H^{\ii} \dfs H_{\kappa_{\iota}^{+}}^{L[E^{P}]}$. By `generic over $H^{\ii}$' we mean that $G^{\ii}_{c}$ intersects every open dense class of $ \mp^{\ii}$ definable over $H^{\ii}$. We note that $\mp^{\ii}$ is itself a proper class of $H^{\ii}$. But $H^{\ii}= L_{\kappa_{\ii}^{+}}[E^{P}]\models
ZF^{-}$, together with a global wellorder of its domain definable over $H^{\ii}$.  Thus given a formula $\psi(v_{0},\vec p)$ with parameters $\vec p \in H_{\ii}$ defining some open dense class $D\sset \mp^{\ii}$, we may define by recursion a maximal antichain $A\sset D$. (6) implies that $|A|\leq \k_\ii$ in $H^{\ii}$ and thus is an element of $H^{\ii}$ by the acceptability of the $\lep$ hierarchy. \qed (12)\\

 We may now define $\Vdash_{\ii} $ the $\mp^{\ii}$-forcing relation 
over $H^{\ii}$. Then we shall have:\\

(13) For $\ii , \nu \in C$, $\ii < \nu$, $\widetilde \pi_{\ii,\nu} : \pa {H^{\ii}, \mp^{\ii}, \Vdash_{\ii}} \longrightarrow_{e} \pa {H^{\nu}, \mpn, \Vdash_{\nu}}$.

\pf By (11) and (12). \qed (13)\\

We let $\langle M_{\infty}, E, \pa{ \pi_{\ii,\infty}}_{\ii\in C} \rangle \dfs \mathrm{dirlim}_{\ii\rightarrow \infty, \ii \in C} $ $ \pa{M_{\ii}, \in, \pa{\pi_{\ii,\nu}}  }_{\ii\leq \nu\in C}$. We may consider $M_\infty$ to be given by an $\in$-relation in the direct limit as some definable (in $V$) class $E \sset V\times V$. This domain we can identify with the  domain $ \pa{H^{\infty},E, \pa{\widetilde \pi_{\ii,\infty}}} \dfs \mathrm{dirlim}_{\ii\rightarrow \infty, \ii \in C} $ $ \pa{H^{\ii}, \in, \pa{\widetilde \pi_{\ii,\nu}} }$, the sole difference being that the maps $\widetilde \pi_{\ii,\nu}$, and so direct limit maps $\widetilde \pi_{\ii,\infty}$ are fully elementary: $\widetilde \pi_{\ii,\infty}:\pa {H^{\ii},\in}\rightarrow_{e} \pa{H^{\infty},E}$. Of course if there were more `ordinals' above $On$ we would say that $\langle H^{\infty},E \rangle$ is isomorphic to a model $\pa{ \widetilde H, \in } \models ZFC^{-} +$ ``$On$ is the largest cardinal''. We define $\mp^{\infty}$ over $\pa{H^\infty,E}$. Note $\mp^{\infty}$ is also a proper class of $H^{\infty}$;
 but nevertheless we can still say that a $\mp^{\infty}$-forcing relation $\Vdash_{\infty} $ for $\pa{H^{\infty},E}$ is definable by taking the direct limit of the relations defined before (13) above. (It would be natural to want to formalise this whole discussion in Kelley-Morse class theory, noting that we have a strong class choice principle in the form of  a global wellorder of $H^{\infty}$ (a model of ``$V=L[E]$'') which is $\pa{H^{\infty}, E}$-definable. Our $KM$-theorem then would additionally talk naturally about all appropriate classes contained in $C$, rather than restricting to $ZFC$-definable ones.)\\

The Mathias condition in this context is obtained by treating $On$ as another indiscernible in $C$:\\

(14) A proper class function $d\in  \prod_{\alpha <\infty}(\mu_{\alpha}\back \bigcup_{\beta < \alpha}\mu_{\beta})^{\omega}$ is $\mp^{\infty}$-generic over $\pa{H^\infty, E}$ if and only if  for every $X\in \prod_{\alpha <\infty}U_{\alpha}$, $X$ coded in $H^{\infty}$, satisfies
 $\bigcup_{\alpha <\infty} (d(\alpha)\back X(\a))$ is finite.\\

Another characterisation of being $\mp^{\infty}$-generic is given below. From now on we let $c = c^{P}$ be the sequence $ \pa{ c^{P}({\a})\mid \a\in On}$ where $ {c^{P}({\a}) }$ is as defined above. Then $\bigcup^{3}c\back \om = P\back P^{\ast}$.

\begin{lemma} {\em (15)} Let $\ii\in C$. Then ${c}\rest \ii$ is $\mp^{P,\ii}$-generic over  $\pa{H^{\ii},\in}$. 
 \end{lemma}
 \pf  The first assertion will follow from the Mathias condition characterized in Theorem \ref{mathias}. But for this we need to observe that
  for every $X\in (\prod_{\alpha <\ii}U_{\alpha})\cap \lep$,
 $\bigcup_{\alpha <\ii} (c(\alpha)\back X(\a))$ is finite. Let $X$ be such a sequence. Then $X\in H^{\ii}$ and as such is in the domain of the direct limit model $M_{\ii}$. We thus have that $X= \pi_{0,\ii}(f_{0})(\kappa_{i^{0}_{0}}, \, \ldots \, , \kappa_{i^{0}_{n(0)}})$ for  some $f_{0}\in M_{0}$, and some indices ${i^{0}_{0}} <\, \cdots \, < {i^{0}_{n(0)}}< {\ii}$. The iteration $\pi_{i^{0}_{n(0)}+1,\ii}: M_{i^{0}_{n(0)}+1} \longrightarrow M_{\ii}
 $ only uses critical points $\kappa_{j} > \kappa_{i^{0}_{n(0)}}$ and first lines up  the next measure onto 
  $\lv{\alpha_{0}}{\omega}$
 where $\alpha_{0}$ is defined to be that least $\alpha$ so that
 $$ \kappa_{{i^{0}_{n(0)}}}
 <  \lv{\alpha}{\omega}
 $$
 and then proceeds with the rest of the iteration to $\ii$. Define $\widetilde X^{\tau}\dfs \pi_{0,\tau}(f_{0})(\kappa_{i^{0}_{0}}, \, \ldots \, , \kappa_{i^{0}_{n(0)}})$ for $\tau >  {i^{0}_{n(0)}}
 $. Then \\
 
 (16) $\pi_{\tau,\ii}(\wt X^{\tau})= \wt X^{\iota}= X
 $.\\
 
 (17) For $\beta \geq \alpha_{0}$ we have: $X(\beta) = \wt X^\tau(\beta)$ for any $\tau \geq \lv{\beta}{ \omega}
 $.\\
 \pf  For such a $\tau$, although $\tau '< \lv{\beta}{ \omega}
 \rightarrow \kappa_{\tau'}< \lv{\beta}{\om}$, (as we are iterating up a smaller measure - meaning not the topmost measure - to $ \lv{\beta}{\om}$) $ \lv{\beta}{\om}$ itself is not a critical point of the iteration, and thus $\kappa_{\tau}>  (\lv{\beta}{\om}^{+})^{\lep}$.  So $\pi
 _{\tau,\ii}\rest L_{\lv{\beta}{\om}^{+}}[E^{P}] = \mathrm{id}\rest L_{\lv{\beta}{\om}^{+}}[E^{P}] $ and so $\wt X^{\tau}(\beta) = \pi
 _{\tau,\ii}(\wt X^{\tau}(\beta)) = \wt X^{\ii}(\beta) = X(\beta).
 $ by (16).
 \qed (17)\\
 
 (18) (i) $c(\alpha_{0})\back  \kappa_{i^{0}_{n(0)}+1}\sset X(\alpha_{0}) $; thus at most finitely many elements of $c(\alpha_{0})$ are not in $X(\alpha_{0})$.
 
 (ii) For $\beta \in (\alpha_{0},\ii)$, $c(\beta)\sset X(\beta)
$.\\
\indent (iii) ${c}\rest [\a_{0},\ii)$ fulfills the condition for $\mp_{\alpha_{0}}^{\ii}$-genericity.

\nod\pf To abbreviate, set $\bar \l =  \lv{\alpha_{0}}{\omega}$ and $j ={i^{0}_{n(0)}+1}$ and $F= E^{M_{j}}_{\nu}
$  the latter the full measure on $\kappa_{j}$ that is being normally iterated up to $\bar \l$. Then $\pi_{j,\bar \l}
(F)=U_{\a_{0}}$ (the full measure on $\bar \l$ in $\lep$ in the notation above). But $\wt X^{j}(\a_{0})\in F$. By normality of the measures in the iteration $\kappa_{j}\in \pi_{j,\bar \l}
(F)$ as well as the intermediate $\kappa_{\tau'}$ for $\tau'\in [j,\bar \l)$. But the latter include, for some $k<\om$, the ordinals $\kappa_{\lv{\a_{0}}{k}}=\lv{\a_{0}}{k} > \kappa_{j}$ which form a co-finite  tail of $c(\alpha_{0})$.\\
For (ii) a similar argument: for $\beta > \a_{0}$, there will be some $j< {\lv{\beta}{1}} (= \kappa_{\lv{\beta}{1}})$ and some $F\in M_{j}
$ such that $\pi_{j,\l_{\om\b+\om}}(F)=U_{\b}$. Setting $F' = \pi_{j,\l_{\om\b+1}}(F)$ this is a full measure on ${\lv{\beta}{1}}$ in $M_{{\lv{\beta}{1}}}$. Then $\wt X^{{\lv{\beta}{1}}}(\beta)\in (\power({\lv{\beta}{1}}))^{M_{{\lv{\beta}{1}}}}$ will have $F'$ measure $1$, and by normality of the iteration from $M_{{\lv{\beta}{1}}} $ to $M_{{\lv{\beta}{\om}}} $, we have for all $k>0$: 
$${\lv{\beta}{k}}\in \pi_{{\lv{\beta}{1}},{\lv{\beta}{\om}}}( \wt X^{{\lv{\beta}{1}}})(\beta) = \wt X^{{\lv{\beta}{\om}}}(\beta) = X(\beta).$$ 
Thus $c(\beta)\sset X(\beta)$.  This concludes (ii) and with (i), (iii) is immediate. \qed (18)\\

 We now repeat the process below $\alpha_{0}$ obtaining a descending chain $\alpha_{0}> \cdots > \alpha_{k}$ of ordinals verifying new, lower, pieces of the form
${c}\rest [\a_{l+1},\a_{l})$ of the condition for $\pa{c({\a})}_{\a<\ii}$. This process will halt with all of $\pa{c({\a})}_{\a<\ii}$ so verified. These details follow.

Then $X\rest [0,\alpha_{0}) \in H^{\lep}_{\lambda_{\om\alpha_{0}}^{+}}$ and as such is in the domain of the direct limit model $M_{\lambda_{\om\alpha_{0}}}$. We thus have that $X\rest [0,\alpha_{0})= \pi_{0,{\lambda_{\om\alpha_{0}}}}(f_{1})(\kappa_{i^{1}_{0}}, \, \ldots \, , \kappa_{i^{1}_{n(1)}})$ for  some $f_{1}\in M_{0}$, and some indices ${i^{1}_{0}} <\, \cdots \, < {i^{1}_{n(1)}}< {\lambda_{\om\alpha_{0}}}$. 
Let $X_{\a_{0}}$ abbreviate $X\rest [0,\alpha_{0})$.

The iteration $\pi_{i^{1}_{n(1)}+1,{\lambda_{\om\alpha_{0}}}}: M_{i^{1}_{n(1)}+1} \longrightarrow M_{{\lambda_{\om\alpha_{0}}}}
 $ only uses critical points $\kappa_{j} > \kappa_{i^{1}_{n(1)}}$ and first lines up  the next measure onto 
  $\lv{\alpha_{1}}{\omega}$ where $\a_{1}$ is the least $\a$ so that $\kappa_{i^{1}_{n(1)}}<  \lv{\alpha}{\omega} $.

Define $\widetilde X^{\tau}_{\a_{0}}\dfs \pi_{0,\tau}(f_{1})(\kappa_{i^{1}_{0}}, \, \ldots \, , \kappa_{i^{1}_{n(1)}})$ for $\tau >  \kappa_{i^{1}_{n(1)}}
 $. Then $\pi_{\tau,{\lambda_{\om\alpha_{0}}}}(\wt X^{\tau}_{\a_{0}})= \wt X^{{\lambda_{\om\alpha_{0}}}}_{\a_{0}}= X_{\a_{0}}
 $. Arguing just as for (17) and (18) above we have:\\

 (19) For $\beta \geq \alpha_{1}$: $X_{\a_{0}}(\beta) = \wt X^\tau_{\a_{0}}(\beta)$ for any $\tau \geq \lv{\beta}{ \omega}
 $.\\

 (20) (i) $c(\alpha_{1})\back  \kappa_{i^{1}_{n(1)}+1}\sset X(\alpha_{1}) $; thus at most finitely many elements of $c(\alpha_{1})$ are not in $X(\alpha_{1})$.
 
 (ii) For $\beta \in (\alpha_{1},\a_{0})$, $c(\beta)\sset X(\beta)
$.\\
\indent (iii) ${c}\rest [\a_{1},\a_{0})$ fulfills the condition for $\mp_{\alpha_{1}}^{\a_{0}}$-genericity.\\

We continue in this fashion defining a descending sequence of critical points $\kappa_{i^{l+1}_{n(l+1)}}< \kappa_{i^{l}_{n(l)}} $, and ordinals $\a_{l+1}<\a_{l}$, and deriving that ${c}\rest [\a_{l+1},\a_{i})$ fulfills the condition for $\mp_{\alpha_{l+1}}^{\a_{l}}$-genericity. We then reach a point where  ${i^{m+1}_{n(m+1)}} =0 $ 
in that for some $f_{m+1} \in M_{0}$ we have that $X\rest [0,\alpha_{m})= \pi_{0,{\lambda_{\om\alpha_{m}}}}(f_{m+1})(\kappa_0)$.
As $\kappa_{0}<\lambda_{0}$ we have  that $c(0)\sset X(0)$ and similarly $c(\b)\sset X(\b)$ for $\beta \in (0,\a_{m})$.  \\

(21)    ${c}\rest \ii$ is $\mp^{\ii}$-generic over $\lep$.\\
\pf 
Setting $\a_{m+1}=0$ we then have: ${c}\rest \ii = \bigcup_{l<m+1} {c}\rest [\a_{l+1},\a_{l}) \cup
{c}\rest [\a_{0},\ii) $, and
$c(\alpha)\sset X(\a)$ for all $\alpha$ not one of the $\a_{l}$.
There are only finitely many $\a_{l}$, so this follows from (18)(i) and $m$ instances of (20)(i). \qed (21)\\

This finishes the Lemma.  \qed (Lemma)\\


\blem (22) If $\ii < \nu \in C$, $\widetilde \pi_{\ii,j}:L[E^{P}]\longrightarrow _{e}L[E^{P}]$, with $\widetilde \pi_{\ii,\nu}(\ii)=\nu=\kappa_{\nu}$ as above arising from the iteration maps $\pi_{\ii,\nu}:M_{\ii}\longrightarrow M_{\nu}$, then with $c\rest\ii$ \etc as above, there is $\overline\pi_{\ii,\nu}\supset \widetilde \pi_{\ii,\nu}$ with 
$\overline\pi_{\ii,\nu}: L[E^{P}][c\rest \ii]\longrightarrow _{e}L[E^{P}][c\rest\nu]$, with 
$\overline\pi_{\ii,\nu}(c\rest \ii) = c\rest \nu $.
\elem 
\pf As $c\rest \ii$ is $\mp^{\ii} $-generic (respectively,  $c\rest \nu$ is $\mp^{\nu} $-generic) over $L[E^{P}]$ and 
$\widetilde \pi_{\ii,\nu}\m``G_{c\rest\ii}\subset G_{c\rest\nu}$, if we define $\overline\pi_{\ii,\nu}(\dot \tau_{G_{c\rest\ii}}) =
\widetilde \pi_{\ii,\nu}(\dot\tau)_{G_{c\rest\nu}}
$, then $\overline\pi_{\ii,\nu}$ will be well-defined and elementary, extending $\widetilde \pi_{\ii,\nu}$. Furthermore
$\overline\pi_{\ii,\nu}(\dot \Gamma^{\ii}_{G_{c\rest\ii}}) =
\overline\pi_{\ii,\nu}(c\rest\ii)= c\rest \nu = \dot \Gamma^{\nu}_{G_{c\rest\nu}}
$.
\qed\\

Consequently:\\

(23) $\bigcup_{\ii\in C}c\rest \ii =c $ is $\mp^{\infty} $-generic over $\pa{H^{\infty},E}$.\\

\pf As we can see, for $c = \bigcup_{\ii\in C} c\rest \ii$ will be $\mp^{\infty}$-generic over the direct limit model $\pa{H^{\infty},E}$, as for any 
$X\in (\prod_{\alpha <\infty}U_{\alpha})\cap \lep$,
 with $X$ coded into $H^{\infty}$ the condition that $\bigcup_{\alpha <\infty} (c(\alpha)\back X(\a)) = \bigcup_{\ii\in C}\bigcup_{\alpha <\ii} (c(\alpha)\back X_ \ii(\a))$  be finite is fulfilled.  \qed (23)\\

We now finish:\\

\pf of Theorem 1.1. For $C$ we take $C_{M_{0}}$ the class of iteration points of the countable mouse $M_{0}$ by its topmost measure.  Given then any cub $P,Q\sset C$ we shall have that there are iteration embeddings $j:L[E_{0}]\longrightarrow L[E^{P}]$
and $k:L[E_{0}]\longrightarrow L[E^{Q}]$. The reals of all such models are thus the same. As the forcings $\mp^{C,\infty}$ (with the obvious definition), $\mp^{{P},\infty}$ and $\mp^{{Q},\infty}$ add no new bounded subsets of their least measurable we shall have that 
$$ L[P]=L[E^{P}][c] \mbox{ and } L[Q]=L[E^{Q}][d]
$$ 
have the same reals, (indeed subsets of $\kappa_{0}$, the least measurable cardinal of $L[E_{0}]$) where $c$, $d$ are  $\mp^{{P},\infty}$- and $\mp^{{Q},\infty}$-generic over  $L[E^{P}]$  , respectively, $L[E^{Q}]$.   By the elementarity of $j,k$ the topmost condition $\mathbbm{1}$ forces the same sentences in the forcing language over the respective models. Hence $ Th(\pa {L[P],\in,P} ) =  Th(\pa {L[Q],\in,Q} )$. 
\qed

\begin{corollary}
If $P$ is appropriate, $L[P]$ is a $\mp^{P,\infty}$-generic extension of its core model - the latter being an iterate of the `minimal' model of a proper class of measurable cardinals, $L[E_{0}]$.
\end{corollary}
\pf With $c$ $\mp^{P,\infty}$-generic over $L[E^{P}]$, $c$ contains none of its limit points. But $P$ is just  $\bigcup^{3}c\back \om$ together with the latter's closure.
It is thus mutually interconstructible with $c$. Hence $L[P]= L[c] =L[E^{P}][c] $. But also $K^{L[P]}=L[E^{P}]$.  \mb \qed\\

With less than a proper class sized $P$ the reader will now see easily that 
similar results apply for set sequences $P,Q\sset C$ of the same limit order type: any two such will have the same reals, the same theories and will look like the same \pr-generic extensions of their inner models $L[E^{P}] $ which now have only a bounded set of measurable cardinals, depending on the length of $P$ or $Q$.

\section{ The H\"artig quantifier model $C(I)$}

We apply the above analysis directly to $C(I)= L[Card]$. However first we show that below $O^{k}$ $C(I)$ computes the canonical inner core model $K$ of $V$.
We then characterise $C(I)$ inside $L[E]$ models. This shows that below $O^{k}$ the H\"artig quantifier picks up all the sets of the model (and this is an equivalence to the non-existence of $O^{k}$).
We shall let $K^{I}= (K)^{C(I)}$.
\begin{theorem}\label{KinC(I)} 
 $\neg O^{k} \Equi K^{I} = K$.
\end{theorem}

\pf $(\Leftarrow)$ follows from the work above: if $O^{k}$ exists, it is not an element of $C(I)$ and hence $K^{I}$  cannot be $K$ which contains $O^{k}$. 
$(\Rightarrow)$: we compare the models $K=M_{0}$ and $K^{I}=N_{0}$ via coiteration  $\cali = \pa{ {M_{\ii}}, \pa{\pi^{M}_{\ii ,j}}_{\ii\leq j\leq\th}, {\kappa_{\ii}} }_{\ii \leq \th}  $ and 
 $\calj = \pa{ {N_{\ii}}, \pa{\pi^{N}_{\ii ,j}}_{\ii\leq j\leq\th}, {\kappa_{\ii}} }_{\ii \leq \th}  $  and with indices $\pa{\nu_{\ii}}_{\ii<\th}$ where in this case ${\th} = \infty$
- the comparison is class length. \\

(i)  $K^{I}$ is  universal.\\
\pf Assume not.  $K^{I}$ can only have boundedly many measurable cardinals. As there can be no truncation on the $N$-side of this coiteration (see \eg \cite{Z02}, Lemma 5.3.1),  there is some least stage $\ii_{0}$ such that for $\mu \geq \ii_{0}$ $\pi^{N}_{\ii_{0},\mu}= id \rest N_{\ii_{0}}$.  That is, all the full measures of $N_{\ii_{0}}$ have been lined up with those of $M_{\ii_{0}}$.     Thereafter $N_{\mu}=N_{\ii_{0}}$. On the $M$-side there may or may not have been a truncation but in any case if $M_{\ii_{0}}$ is still a proper class, there is an initial segment  of $M_{0}$, some $M_{0}| \tau$ say so that the coiteration of $M_{0}|\tau$ with $N_{0}$ yields the same outcome with the same indices and ultrapowers taken on both sides. We may this assume that $M_{0}$ is replaced by such a $M_{0}|\tau$. The point is that some set sized mouse will eventually iterate using repeatedly only 
some  filter $F_{{\ii}} = E^{M_\ii}_{\n_{\ii}}$ and its images, with critical point $\k_{\ii}$, for $\ii \geq \ii_{0}$, past $N_{\ii_{0}}$ leaving this model as $N_{\infty}$ behind. Let  $
\pa{\l_{\a}\mid \a < \omega_{1}}$ increasingly 
enumerate the next $\om_{1}$ $V$-cardinals above
$|M_{\ii_{0}}|^{+ }=\l_{0}$, 
and let their supremum be $\lambda$. 

Then (a) the sequence $\pa{\l_{\a}}_{\a<\om_{1}}\in C(I)$; (b)  as the cardinality of $M_{\ii_{0}}<\lambda_{0}$ each of the $\l_{\a}$  satisfy that $\k_{\l_{\a}} = \l_{\a}$ and $\pi^{M}_{{\ii_{0}},\l_{\a}}(\k_{\ii_{0}})= \k_{\l_{\a}}= \l_{\a}$. Consequently the filter $F_{\l}= E^{M_\l}_{\n_{\l}}$ to be used at stage $\l$ is generated by the final segment filter using the sequence $\pa{\k_{\ii}}_{\ii_{0}<\ii<\l}$, but also by the subsequence $\pa{\k_{\l_{\a}}}_{\a<\om_{1}}= \pa{{\l_{\a}}}_{\a<\om_{1}}$. But at this stage $\l$ we have $(\power(\l))^{M_{\l}}= (\power(\l))^{N_{\ii_{0}}}$
We further have (c): the cardinals $\lambda_{\a}$ are all fixed points of the embedding $\pi^{N}_{0,\ii_{0}}$. We may thus,  in $C(I)$, define $\bar F$ on $(\power(\l))^{N_0}$ using the same final sequence $\pa{\k_{\l_{\a}}}_{\a<\om_{1}}$. Thus $X\in \bar F \equi \pi^{N}_{0,\ii_{0}}(X)\in F_{\l}$.
This is an $N_{0}$-normal amenable measure on $\lambda$, which is again $\om$-complete. We have a contradiction as on the one hand $K^{I}$ is universal in $C(I)$ (it is the actual core model of $C(I)$), and thus by the theory of such models all  $\omega$-complete normal measures amenable to it are on its $E^{K^I}$ sequence; whilst on the other all the measurable cardinals of $K^{I}$ are strictly below $\k_{\ii_{0}}$ which is less than $\l$. \qed (i)\\


(ii) $K= K^{I}$.

\pf As we are below $O^{pistol}$ (the sharp for  the least model of a strong cardinal) it is a theorem of Jensen (see \eg \cite{Z02}, Thm.7.4.9) that any universal weasel $W$, and by (i) $K^{I}$ is such, is a simple iterate of $K$ in which $W$ does not move. If $K^{I}$ contains no measurable cardinals then the result is proven: there are no measurables in $K$ to iterate. 

Suppose $K\neq K^{I}$ for a contradiction. In the comparison of $K$ with $K^{I}$ let the first measurable to be moved on the $K$-side be $\kappa$, and let us suppose it to be iterated up to the measurable cardinal $\k^{I}$ in $K^{I}$. Suppose there is a further measurable cardinal in $K$ above $\k$ which has critical point $\l$ (where we take $\l$ least). Then the measure on $\l$ here is to be iterated up to some measurable $\l^{I}> \k^{I}$ in $K^{I}$.
(The case of only the one $\k$ measurable in $K$ will be left to the reader.)
So we suppose $K$ (and so $K^{I}$) has at least two measurable cardinals.

Note that $\k \leq \k^{I}$ and $\l \leq \l^{I}$. Let $\l^{I}<\mu$
where the latter is a strong limit $V$-cardinal of $C(I)$-cofinality greater than $\tau$, where we set $\tau= |\k|^{+}$. 

In $K$ iterate the measure on $\l$ $\mu$-times, up to $\mu$ and do the same in $K^{I}$ sending $\l^{I}$ to $\mu$. Let the resulting models be $\bar K$ and $\bar K^{I}$ on the respective sides. Further let $M= \bar K_{\mu}$ and $N=\bar K^{I}_{\mu}$ be the initial segments of the two models cut down to $\mu$.   These are $ZFC$-models. To compare these two all we have to do is iterate the single measure $F_{\k}$ on $\k$ in $M$  up to $F^{I}_{\k^{I}}$ on $\k^{I}$ in $N$. Let $\sigma: M\longrightarrow N$ be this iteration map. If $\la \lambda_{\a}\mid \a < \tau\ra$ strictly increasingly enumerates the next $\tau$ many successor elements of $Card$ above $\k^{I ++}$ then all such $\l_{\a}$ are less than $\mu$ but are fixed points of the iteration map $\sigma$. Fix a set of definable skolem functions for $N$, and so for $M$ too by elementarity. Let $H\prec N$ be the skolem hull of: $\kappa \cup \{ \l_{\a}\}_{\a < \tau}\cup {F^{I}_{\k^{I}}}$. Then $H\cong H'$ where the latter is  
the skolem hull of: $\kappa \cup \{ \l_{\a}\}_{\a < \tau}\cup {F_{\k}}$ in $M$. (Note that we have $\sigma(F_{\kappa}) = F^{I}_{\k^{I}}$ and $H\cap (\kappa, \kappa^{I})=\emp$ (the latter as $\sigma(\k)=\k^{I}$).) Furthermore $H$ is definable in $C(I)$ since $N$ and those components are. Let $\pi : H \longrightarrow P$ be the Mostowski-Shepherdson Collapse, and
thus $|P| \geq \tau$. 

 Check that $\pi(F^{I}_{\k^{I}})$ collapses to $F_{\kappa}\sset P$, all inside $C(I)$. Hence, as we have $E^{K}\rest \k^{+} = E^{K^{I}}\rest \k^{+}$ we have $Q= \pa{J^{E^{K}}_{\k^{+}}, E^{K}\rest \k^{+}, F_{\kappa}} \in C(I)$, and we may then proceed to build, in $C(I)$,  a universal class model $W$ with this structure as an initial segment. (In Jensen's nomenclature $Q$ is a `strong mouse', \cf \cite{Z02} Section 7.1 and Lemma 7.1.1).  This is a contradiction since in $C(I)$ $K^{I}$ must simply iterate up to $W$ and thus no such $W$ can have a measurable cardinal on an ordinal less than $\k^{I}$.\\ \qed ((ii) and Theorem)\\

We then have easily that inside canonical models, if they are not too large then they are their own H\"artig quantifier models:

\begin{corollary}\label{CorKinC(I)} $(V=L[E])$ \,\, $\neg O^{k}\equi V=C(I)$.
\end{corollary}
\pf Again if $O^{k}$ exists in $V$ we have ensured it is outside of $L[Card]= C(I)$. If $V=L[E]$ then $V=K$. If additionally $\neg O^{k}$ then we also have $K\sset C(I)$. \qed\\

For an inner model $W$ let $WCL(W)$ mean that for $\all \tau
(\tau \mbox{ a singular cardinal } \Imp (\tau^+ = (\tau^{+})^{W})$. 

\begin{corollary}  $\neg O^{k}\equi WCL(C(I))$.
\end{corollary}
\pf This is immediate since $WCL(K)$, and then $WCL(C(I))$ will hold if $K\sset C(I)$. \qed

\begin{corollary}\label{Cor5.4} Assume $ O^{k}$ exists. Then $C(I)$ is a $\mp^{Card,\infty}$-generic extension of its core model $K^{I}$, where $\mp^{Card,\infty}$ is as defined above as  $ \mp^{P,\infty}$ for $P=Card$. 
\end{corollary}

\begin{corollary} Assume $ O^{k}$ exists. Then $C(I)\neq HOD$.
\end{corollary}
\pf $ O^{k}\notin C(I)$, whilst $ O^{k}\in HOD$. \qed 

\begin{corollary} Assume $ O^{k}$ exists.  $C(I)^{C(I)}= K^{C(I)}$.
Consequently ${C(I)^{C(I)}} \models \m`` V= {C(I)}\mq.$
\end{corollary}
\pf Note that by Cor. \ref{Cor5.4} $C(I)$ has the same cardinals as its core model $K^{I}$.   Also $K^{I}$ satisfies ``$V=L[E]$'', thus by Cor. \ref{CorKinC(I)}, we have $C(I)^{C(I)} =C(I)^{K^I}= K^{I}$.
\qed\\

The latter may consistently  fail if $\neg O^{k}$: let $M$ be the forcing extension of $L$ that adds a Cohen real $r$, and then collapses $\aleph_{2n+1}$ to $\aleph_{2n}$ iff $n\in r$. Then $(V=C(I)= C(I)^{C(I)})^{M}$. But $K^{M}=L$. 

\section{The $\tmop{Cof}_{\omega}$ model $\cs$}

We briefly make a few comments on the relationship between the H\"artig quantifier model $C(I)$ and the $\tmop{Cof}_{\omega}$ model $\cs$ of \cite{KeMaVa2016}. For our purposes here we let $\tmop{Cof}_{\omega} = \{\a\mid \tmop{cf}(\a)=\om \}$, and then $\cs = L[\tmop{Cof}_{\omega}]$. We show these models differ in that $O^{k}\in \cs$ (if it exists) whilst we have shown this must fail for $C(I)$. 
We first note as an aside a generalisation of an argument of \cite{KeMaVa2016} from a single measure to  a sequence of such.

\begin{theorem} Assume $V= L[E]$ has measurable cardinals $\pa{\k_{\ii}\mid \ii < \theta < \k_{0}}$ with measures $E_{\k_{\ii}}$.  
 Let $L[E']$ be the simple iteration of $L[E]$ where each measure is iterated in turn $\om^{2}$ times with iteration points $\pa{\k_{\ii}^{\a}\mid \a < \om^{2}}$. Then 
${\cs} = L[E'][\pa{\pa{{\k_{\ii}^{\om\cdot n}\mid 0< n < \om}}\mid \ii < \th}]$.

\end{theorem}
\pf  The assumption 
 on the length of the sequence of ensures that all the measurable cardinals $\k_{\ii}$ are discrete ordinals: there are no measurable limits of measurables, and thus those arguments in  \cite{KeMaVa2016} can be straightforwardly deployed for each cardinal in turn.
 \qed \\

We first show by methods of Theorem \ref{KinC(I)} that below $\ok$ $\ks$ is universal. We then show that if $\ok$ exists it must be in $\ks$.  However first we give two lemmata about cofinalities of regular cardinals in iterates of mice.
We do the  ``$n=0$'' example of $M_{0}$ the $\ok$ mouse in detail first, and just state the generalisation for the case $n>0$ afterwards. 
\blem Let $\pi_{0,\th}:M_{0}\longrightarrow M_{\th}$ be a simple normal iteration of $M_{0}$, the $\ok$ mouse, (\ie without any truncations), with critical points $\pa{\k_{\ii}\mid \ii < \th}$.  Suppose that:\\ (i) $M_{\th}\models \mbq \k \mbox{ is inaccessible, but not measurable''}$ ; (ii) $\k \neq \k_{\ii}$ for any $\ii<\th$.\\ Then $cf^{V}(\k)=\om$.

\elem
\pf By induction on $\th$. Suppose $\th = \th_{0}+1$. 
If $\k < \k_{\th_{0}}$ then the result follows from the inductive hypothesis as $\pi_{\th_{0},\th_{0}+1}(\k)=\k$. By assumption $\k = \k_{\th_{0}}$ is ruled out. For $\k > \k_{\th_{0}}$ we use the following observations:

{\em Claim (1)} If $M_{\ii}= \pa{J_{\a_{\ii}}^{E^{M_{\ii}}},E^{M_{\ii}}, F_{\ii}}$ then $cf(\a_{\ii})=\om$. Moreover any $\l\in (RegCard)^{M_{\ii}}$ with $\l > \k_{\ii}$ has $V$-cofinality $\om$.

\pf We claim there is a $\Sigma_{1}^{M_{\ii}}$-definable sequence $\pa{\a^{n}_{\ii}}_{n<\om}$ which is increasing and cofinal in $\a_{\ii}$. Let $\a^{0}_{\ii}$ be the least $\a$ greater than $\ki$ so that $\ex \beta \in (\ki,\a ) \, J_{\beta}^{E^{\mi}} \models ZF^{-}\wedge F_{\ii}\cap J_{\beta}^{E^{\mi}} \in J_{\a}^{E^{\mi}} $. By amenability of $\pa{M_{\ii},F_{\ii}}$ $\a^{0}_{\ii}$ and each $\a^{n+1}_{\ii}$ to follow is well defined. 
For this, let $\a^{n+1}_{\ii}$ be the least $\a$ greater than $\a^{n}_{\ii}$ so that 
$F_{\ii}\cap J_{\a^{n}_{\ii}}^{E^{\mi}} \in J_{\a^{n+1}_{\ii}}^{E^{\mi}} $.  Then $\sup_{n}\a_{\ii}^{n} = \a_{\ii}$: for if not then setting $\a' =\sup_{n}\a_{\ii}^{n} $ we should have that 
$M'= \pa{J_{\a'}^{E^{M_{\ii}}},E^{M_{\ii}}, F_{\ii}\cap J_{\a'_{\ii}}^{E^{\mi}}}$ is an iterable premouse with an amenable topmost measure illustrating that it is a measurable limit of measurables. $M'$ is thus, being an initial segment of $M_{0}$,  in the mouse ordering $<^{\ast}$ below $M_{0}$ which was defined to be the least such mouse. Contradiction! We thus have the first sentence of the Claim. $M_{\theta}$, being a simple iterate of $M_{0}$, has the first $\Sigma_{1}$-projectum dropping below its topmost critical point to $\omega$. The same is true of $M_{\ii}$ and thus the latter is the $\Sigma_{1}$-Hull of $\ki$ in $M_{\ii}$ (we write this as $M_{\ii}\dfs \Sigma_{1}\mbox{-SH}^{M_{\ii}}(\ki)$). Let $H^{\ii}_{n}\dfs \Sigma_{1}\mbox{-SH}^{M_{\ii}|\a^{n}_{\ii}}(\ki)$, where
$M_{\ii}|\a^{n}_{\ii} =  \pa{J_{\a^{n}_{\ii}}^{E^{M_{\ii}}},E^{M_{\ii}}, F_{\ii}\cap J_{\a^{n}_{\ii}}^{E^{\mi}}}$. Then $H^{\ii}_{n}\in M_{\ii}$ and thus $\tau^\ii_ n( \l) \dfs H^{\ii}_{n} \cap \lambda$ must be bounded in $\l$ for any $M_{\ii}$-regular $\l > \ki$. However $\sup _{n} \tau^\ii_ n( \l) =\l$ and the Claim then is proven.\\ \qed (Claim (1))\\

The Lemma follows then in the successor case, since if $\k>\k_{\th_{0}}$ then as the iteration is  normal $\k \neq \k_{\ii}$ for any $\ii < \th_{0}$. Suppose now $\tmop{Lim}(\th)$. Let $\k$ satisfy (i) and (ii) of the Lemma. As $M_{\th}$ is a direct limit model, for an $\ii_{0}<\th$ let $\l_{i}\in M_{\ii}$ be such that $\pi_{\ii,\th}(\l_{\ii})=\k$, for $\ii\in (\ii_{0},\th)$. By elementarity, for such $\ii$, $(\l_{\ii} \mbox{ is not measurable but is inaccessible})^{M_{\ii}}$. Consequently $\l_{\ii}\neq \k_{\ii}$. If for some such $\ii$ $\l_{\ii}<\k_{\ii}$ then the conclusion of the lemma holds by the inductive hypothesis in $M_{\ii}$ (as $\pi_{\ii,\th}(\l_{\ii})=\l_{\ii}<\k_{\ii}$). So we may assume $\l_{\ii}>\k_{\ii}$.  However now we may form, as in the proof of the claim above, the ordinals $\tau^\ii_ n( \l_{\ii}) $. The $\Sigma_{1}^{M_{\ii}}$-definability of the sequence $\pa{\a^{n}_{\ii}}_{n}$ yields the same for the sequence $\pa{\tau^\ii_ n( \l_{\ii})}_{n} $. Although the whole iteration $\pi_{\ii,\th}:M_{\ii}\longrightarrow M_{\th}$ is not internally definable in $\mi$, each ultrapower stage $\pi_{\ii,\ii + 1}$ by the measure $E\in M_{\ii}$ with critical point $\ki$, has $\pi_{\ii,\ii + 1}\mbq \l_{\ii}$ cofinal in $\pi_{\ii,\ii + 1}(\l_{\ii})=\l_{\ii+1}$. Further $\pi_{\ii,\ii + 1}(\tau^\ii_ n( \l_{\ii}))=\tau^{\ii+1}_ n( \l_{\ii+1})$. Proceeding to the direct limit we see that the image of this $\om$-sequence will be an $\omega$-sequence cofinal in $\k$. This concludes the limit case and the lemma. \qed (Lemma)\\

We state here the generalisation of this for other mice in this region. We say that an iteration $\sigma: P_{0}\longrightarrow P_{\theta}$ has ``no drops'' if there are no truncations in the iteration, and there are no ``drops in degree'', \ie if $n<\om$ is such that $\rho^{n+1}_{P_{0}} \leq \k_{0}< \rho^{n}_{P_{0}}$ then $\all \ii<\theta (\rho^{n+1}_{P_{i}} \leq \k_{\ii}< \rho^{n}_{P_{i}})$. (Hence the level at which the fine-structural ultrapowers are taken remains constant.) The lemma is proven by similar reasoning to the previous one, which is only an instance of the next with $n=0$.

\blem Let $P \leq^{\ast}M_{0}$. Let $\sigma: P_{0}\longrightarrow P_{\theta}$ have no drops in the above sense. Let $n$ be least with $\rho^{n+1}_{P_{0}} \leq \k_{0}$.
Let $\tau = \tmop{cf}^{V}(\rho^{n}_{P_{0}})$. Suppose 
$\rho^{n+1}_{P_{\ii}} < \k < \rho^{n}_{P_{\ii}} $
and that:\\
  (i)\, $P_{\theta}\models \m``\k$ is inaccessible but not measurable''; (ii)\, $\k \neq \k_{\ii}$ for any $\ii<\th$.\\
   Then $cf^{V}(\k)=\tau$.
\elem

\pf 
We just sketch the main point: although $n$ may be non-zero, the iteration map $\pi_{\ii,\ii+1} : M_{\ii} \longrightarrow M_{\ii+1}$ at each stage is $\Sigma^{(n)}_{0}$-preserving and cofinal at the $n$'th projectum level (and so $\Sigma^{(n)}_{1}$-preserving) (see \cite{Z02}). The map restricted to $H^{M_{\ii}}_{\rho^{n}_{M_{\ii}}}$ is thus cofinal into $H^{M_{\ii+1}}_{\rho^{n}_{M_{\ii+1}}}$. (We recall that  $\rho^{n+1}_{M_{\ii}}$ equals $\rho^{n+1}_{M_{0}}$ for any $0< \ii\leq \theta$.)
 Moreover we may pick a definition for a $\Sigma^{(n)}_{1}$-definable  partial, but cofinal, map $\gamma:\rho^{n+1}_{M_{{0}}}\longrightarrow \rho^{n}_{M_{{0}}}$ which thus furnishes an increasing sequence $\g^{0}_{\eta}$ for $\eta < \tmop{o.t.}(\ran(\gamma)) $
cofinal in $
\rho^{n}_{M_{0}}$. 
The range of $\gamma$ is thus preserved by this definition throughout the iteration as a $\Sigma^{(n)}_{1}$-definable set which we may write in increasing order as $\g^{\ii}_{\eta}$ for $\eta < \tmop{o.t.}(\ran(\gamma)) $, with $\sup_{\eta }\g^{\ii}_{\eta} = \rho^{n}_{M_{\ii}}$ and $\pi_{0,\ii}(\g^{0}_{\eta}) = \g^{\ii}_{\eta}$. It is clear then that $cf^{V}(\rho^{n}_{M_{0}}) = cf^{V}(\rho^{n}_{M_{\ii}})$ for all $\ii \leq \theta$.
Now we can finish off as in the last lemma defining for any regular $\kappa$ in the interval, sequences cofinal in the  ordinal $\kappa$ by defining appropriate skolem hulls in the successor case, and the mechanism of $\g^{\ii}_{\eta}(\lambda ^{\ii})$ as analogues of the $\tau^{\ii}_{n}(\lambda)$  for the direct limit argument \etc
\qed\\

\begin{theorem}\label{KinC*} 
 $\neg O^{k} \rightarrow K^{\ast} \dfs (K)^{C^{\ast}}$ is universal; thus $K^{\ast}$ is a simple iterate of $K$.
\end{theorem}

\pf We argue as in the proof of (i) of Theorem \ref{KinC(I)} assuming that $K^\ast$ is not universal for a contradiction, and thus it only has boundedly many measurable cardinals again. We additionally require that $\io$ is such that there are no further drops on the $M$-side for $\ii \geq \io$. Let $\delta \dfs  \tmop{cf}^{V}(\k_{\io}^{+})^{M_{\io}}$. Further, instead of setting $\pa{\lambda_{\a}}_{{0<\a<\om_{1}}}$ as the next $\om_{1}$ $V$-cardinals above $|M_{\io}|^{+}$, we take them as the next ordinals in increasing order satisfying:

\nod (a) $\tmop{cf}^V(\l_{\a}) = \omega, \mbox{ if } \delta \neq \omega $; or $\neq \om $ if $\delta = \om$; \\
(b) $(\l_{\a}$ is inaccessible$)^{\ks}$.\\

\nod {\em Claim (i)} 
Each $\l_{\a}$ is a fixed point of $\pi^{N}_{0,\io}$.\\
\pf As each of the $\lambda_{\alpha}$ are inaccessible in $N_{0}=\ks$ they would be trivially fixed points for any iteration of $N_{0}$ by measurable cardinals below $\lambda_{0}$ if the iteration were to be internally definable in $N_{0}$. But seemingly there is no guarantee of this. So instead we deploy the {\em universal iteration} of Definition \ref{defuniv}. Let $\theta$ be a regular cardinal of $\ks$ bounding the measurable cardinals of $N_{\ii_{0}}$ and so of $N_{0}= \ks$. By increasing the choice of $\lambda_{0}$ if necessary, we shall assume without loss of generality that $\l_{0} > \theta$. Let $\gamma$ be in $\ks$ an inaccessible limit of inaccessibles. Fix then a universal iteration of length $\theta+1$ as defined with starting model  $Q_{0}= \ks\rest\gamma$. Although the iteration $\pi^{N}_{0,\io}:N_{0}\longrightarrow N_{\io}$ is not defined in $\ks$, Lemma \ref{fixed1.5} then shows that the $\lambda_{\a}$ are all fixed points of this map.\\


 {\em Claim (ii)} Each $\l_{\a}= \k_{\beta(\a)}$ for some critical point in the iteration $\pi^{M}_{\io,\infty}$ of $M_{\io}$ by the top measure $F_{\io}$ on $\k_{\io}$.\\
\pf  The instances of (ii) follow from the last two lemmata above.\qed\\

 We then finish off as follows: let $\l \dfs \sup\{\l_{\a}\}_{\a < \om_{1}}$ then $\l = \k_{\l}$ and $F_{\l}$ is generated by the final segment filter on $\pa{\l_{\a}}$. As these are fixed points of the embedding $\pi^{N}_{0,\io}$ we can define as before, in \cs, $\widetilde F$ as this final segment filter generated by $\pa{\l_{\a}}$ on $\power(\l)^{N_{0}}$. As $\widetilde F$ is then an $\om$-complete measure, we get a contradiction as before.


The next corollary  is just a particular example of the above.
\begin{corollary}
If $\neg O^{\dagger}$ but there is an inner model $L[U]$ (say with $U$'s critical point on the least possible ordinal), then $\ks$ is an iterate of $L[U]$.
\end{corollary}

\begin{corollary}
If there is an inner model with a proper class of measurable cardinals, then there is such an inner model in \cs.
\end{corollary}

Then the following ensures that \cs \, must be different from $C(I)$.

\begin{theorem} {If $O^{k}$ exists, then it is in $C^{\ast}$.}

\end{theorem}

\pf
Assume for a contradiction that $\ok\notin \cs$. We coiterate $P_{0}\dfs K$ with $N_{0} \dfs \ks$ to models $(P_{\infty}, N_{\infty})$. By assumption $\ks \leq^{\ast} L[E_{0}]$, the latter again  the model left behind by the iteration out of $M_{0}$'s top measure.\\
{\em Case 1 $\ks$ has a proper class of measurable cardinals. }\\As $\ok$ exists it is in $K$ and indeed appears as an initial segment of $K$ on the $E^{K}$ sequence. The  coiteration immediately starts with a truncation to a $P_{0}^{\ast} = M_{0}$ of $P_{0}$,  followed by an ultrapower $\pi^{P}_{0,1}:M_{0}\longrightarrow P_{1} $ and thereafter we have a comparison of $P_{1}$ with $N_{0}$.  Thus $\pi^{P}_{1, \infty}:M_{0}\longrightarrow P_\infty $ is a simple normal iteration of $M_{0}$ that generates $N_{\infty}$. In this iteration $\ks=N_{0}$ does not move, and thus $N_\infty = N_{0}$, by Lemma \ref{minimality}.
Now consider $\vec c $ an increasing $\om$-sequence of ordinals $\nu_{n}$ that (i) have uncountable $V$-cofinality; (ii) are limits of measurable cardinals in $K^{\ast}$, and (iii) are inaccessible in $K^{\ast}$. Such a sequence must exist
as there is a cub class of $\ii$ where the topmost measure $F_{\ii}$ of $P_{\ii}$ is used to form an ultrapower at stage $\ii$ (and this leaves behind a non-measurable but inaccessible limit of measurables in $L[E^{K^{\ast}}]$). 
Such can be found in $\cs$ as $C^{\ast}= L[\tmop{Cof}_{\om}]$. But conversely any $\nu_{n}$ satisfying (i)-(iii)  must itself be a critical point $\k_{\ii(n)}$, where by (i) and (ii) of the Lemma 6.2 the step $\pi_{\ii(n),\ii(n) + 1}$  has to be  an ultrapower step by the topmost measure $F_{\ii(n)}$ of $P_{\ii(n)}$. If ${\ii^\ast} = \tmop{sup}_{n}\{{\ii(n)}\}$, 
then in the direct limit model
the topmost measure $F_{\ii^{\ast}}$ of $ P_{\ii^{\ast}}$ on $\power(\k_{\ii^{\ast}})\cap K^{\ast}=\power(\k_{\ii^{\ast}})\cap P_{\ii^{\ast}}$ is generated by the final segment filter on $ \pa{\k_{\ii(n)}}_{n}\in \cs$. But $F_{\ii^{\ast}}$ is then in $\cs$, and so $P_{\ii^{\ast}} \in \cs $. But $M_{0}= O^{k}= core(P_{\ii^{\ast}})$, that is, it is the (transitive collapse of) the $\Sigma_{1}$-SH$^{P_{\ii^{\ast}}}(\emp)$, and thus is also in $\cs$.

\nod{\em Case 2 Otherwise} \\
We argue that this case cannot occur. If it did, then $\ks$ has a bounded set of measurable cardinals at most. Now argue as in the proof of Theorem \ref{KinC*}. For some $\io$ there are no further truncations on the $P$-side of the iteration for $\ii\geq \io$. We take 
$\l_{\a}$ for $0\leq\a<\om_{1}$ an ascending sequence in $\cs$ of ordinals satisfying (a) and (b) there. 
(Again, apply the arguments using the universal iteration and Lemma \ref{fixed1.5} that $\pi^{N}_{0,\eta}(\l_{\a})=\l_{\a}$.)
We again define an $\widetilde F$, an $\om$-complete measure on $\power(\l)^{\ks}$ for $\l = \sup_{\a}\l_{\a}$, with $\widetilde F \in \cs$. This is a contradiction just as before. So Case 2 cannot occur. \qed\\





As discussed above, the exact nature of $\cs$ remains open, but the above methods illustrate starkly how they do not apply to the least sword mouse $O^{s}$.



\end{document}

%% file: pdwstyle.tex
\newcommand{\cali}{\ensuremath{\mathcal{I}} }
\newcommand{\calj}{{\mathcal{J}}}

\newcommand{\call}{ {\mathcal{L}}}

\newcommand{\phiabp}[2]{\ensuremath{\Phi_{\beta +1}^{(  \alpha )}}}

\newcommand{\nod}{{\noindent}}

\newcommand{\cs}{\ensuremath{\ul \sigma \ur}}

\newcommand{\co}[1]{\ensuremath{\ulcorner #1 \urcorner}}

\newcommand{\bu}{\ensuremath{\bullet}}

\newcommand{\rest}{{\upharpoonright}}
\newcommand{\emp}{{\varnothing}}
\newcommand{\pa}[1]{\ensuremath{\langle #1 \rangle}}

\newcommand{\calu}{\ensuremath{\mathcal{U}}}

\newcommand{\om}[1]{\ensuremath{\omega_{1} ( #1 ) }}

\newcommand{\la}{{\langle}}
\newcommand{\ra}{{\rangle}}
\newcommand{\ran}{\text{ran}}
\newcommand{\dom}{\text{dom}}
\newcommand{\back}{{\backslash}}
\newcommand{\power}{P}

\newcommand{\ie}{{\itshape{i.e.}}{\hspace{0.25em}}}
\newcommand{\eg}{{\itshape{e.g.}}{\hspace{0.25em}}}
\newcommand{\etc}{{\itshape{etc.}}{\hspace{0.25em}}}

\newcommand{\cf}{{\itshape{cf. }}{\hspace{0.25em}}}
\newcommand{\mm}{{\itshape{m.m.}}{\hspace{0.25em}}}

\newcommand{\all}{{\forall}}
\newcommand{\Equi}{{\,\Longleftrightarrow\,}}
\newcommand{\equi}{{\,\longleftrightarrow\,}}
\newcommand{\ex}{{\exists}}
\newcommand{\Imp}{{\,\Rightarrow\,}}
\newcommand{\sset}{ { \subseteq } }

\newcommand{\wt}{\tilde{\}}}

\newcommand{\fin}[1]{\ensuremath{[  ]^{< \omega}}}

\newcommand{\qed}{\hfill Q.E.D.}

\newcommand{\pf}{{\noindent}{\textbf{Proof: }}}
\newcommand{\dfs}{\ensuremath{=_{\ensuremath{\operatorname{df}}}}}

\newcommand{\pr}{\ensuremath{\prec}}

\newcommand{\oddpagetext}[1]{\newcommand{\pageoddheader}{{\small }}}
\newcommand{\evenpagetext}[1]{\newcommand{\pageevenheader}{{\small }}}

